# G-KdVNet: ANN–ADM Surrogate for Geophysical KdV Equation

[1]Mrutyunjaya Sahoo, [2]Arup Kumar Sahoo*, [3]Snehashish Chakraverty,

[1]Centre for Data Science, Department of Computer Science and Engineering, Siksha 'O' Anusandhan Deemed to be University, India.

[2]Department of Computer Science and Engineering, Siksha 'O' Anusandhan Deemed to be University, India.

[1,3] Department of Mathematics, National Institute of Technology Rourkela, India.

**Emails**: mrutyunjayasahoo@soa.ac.in[1], arupkumarsahoo@soa.ac.in[2] ,sne_chak@yahoo.com3

*C.A: arupkumarsahoo@soa.ac.in

## Abstract

This research examines the influence of the Coriolis parameter on the behaviour of the geophysical Korteweg-de Vries (KdV) equation. To efficiently approximate its solution, a novel surrogate framework, termed *G-KdVNet*, is proposed by integrating artificial neural networks with the Adomian decomposition method (ADM). In the proposed approach, ADM is first employed to generate reliable semi-analytical solution data, which are subsequently used to train the neural network model. The developed model demonstrates strong predictive capability in capturing the nonlinear dynamics of the KdV system. Numerical results indicate that the proposed model achieves improved accuracy compared with conventional baseline methods, with absolute errors of the order of $10^{-3}$for unseen data. The results suggest that the proposed ANN–ADM surrogate offers an efficient and accurate alternative for solving nonlinear geophysical models, with potential applicability to a broader class of dispersive wave equations.

## 1. Introduction

Many physical and engineering problems are modeled using partial differential equations (PDEs). To gain deeper insights into real-world phenomena, the pursuit of both exact and numerical solutions of PDE has become indispensable, and these have been successfully extracted using a variety of methods. Some examples include the finite element method, finite volume method, homotopy perturbation method, homotopy transformation method, Adomian decomposition method, etc. Fluid dynamics is an emerging field of engineering whose governing equations can be modeled by PDEs.

Artificial intelligence (AI) has shown its potential to revolutionize many areas, including Computational Fluid Dynamics (CFD). Some of the AI techniques, including Artificial Neural Networks (ANN), Robotics, Natural Language Processing, and Fuzzy logic, have shown exponential growth in the last decade by surpassing human levels of accuracy [1]. AI-based numerical computing paradigms have been frequently used to study the behavior of dynamic problems in the last few years. As such, a large amount of experimental data is generated. To effectively manage and analyze this data, data-driven models, particularly ANN models, have emerged as indispensable tools. So, there is a need for time to explore and exploit the numerical technique based on intelligent computing paradigms to solve and analyze fluid problems.

The Korteweg-de Vries (KdV) equation is useful for the evolution of small-amplitude and long-wavelength shallow water. When waves propagate in shallow water, they are affected by the water depth and bottom topography, which can cause the wave speed to vary with depth. As a result, the wave profile changes as it propagates, leading to the phenomena such as wave breaking and wave dispersion. To understand these phenomena, the KdV equation was used to describe the behavior of long waves in shallow water. From the solution of the KdV equation, one may predict the behavior of waves over time and space, including the formation of solitons. Here, we have taken an expanded variant of the KdV equation that incorporates the Coriolis term, known as the Geophysical Korteweg-de Vries (GKdV) equation. This is more significant for ocean waves on a global scale.

The remainder of the paper has the following format: In Section 2, a brief literature study is presented. For the sake of completeness, some preliminaries related to the ANN and GKdV equations are discussed in Section 3. In Section 4, a brief description of the proposed methodologies and their implementation has been presented. Some numerical results and the effect of the addition of the Coriolis constant to the governing equations are included in Section 5. In Section 6, a summary of the findings from the current investigation is presented.

## 2. Related Studies

The KdV equation was first presented in 1895 by Korteweg and De Vries [2]. Various research works have been published on KdV-like equations over the past few years. Subsequently, scientists made different modifications to this equation. In this regard, Johnson and Holm [3] conducted a study on the Camassa-Holm equation, the KdV equation, and other related models for water waves. Kudryashov [4] provides a traveling-wave general solution to the KdV and KdV-Burger equations. Seadawy et al. [5] have used an algebraic method to observe the traveling wave solution to different kinds of KdV equations. To solve time-dependent Schrodinger-KdV equations, Cai et al. [6] used the composition method to develop temporal second and fourth-order schemes. Karunakar and Chakraverty [7] have modeled a variety of wave phenomena using the homotopy perturbation transform method (HPM), specifically for two distinct time-fractional KdV equations of orders three and five, respectively. Aljahdaly et al.

[8] used the direct algebraic method to find an analytic solution for the generalized seventh-order KdV equation that occurs in shallow water waves (SWWs). Additionally, Sahoo and Chakraverty [9], [10] have employed a hybrid method, specifically the Sawi transform-based HPM and only HPM, to obtain approximate solutions for one-dimensional coupled equations of SWWs in the context of tsunami predictions.

ADM is one of the most important methods for dealing with these nonlinear models. In the early 1990s, Adomian suggested this novel method for solving some functional equations [11], [12]. There are significant benefits of using ADM instead of standard numerical approaches. The approximation has the advantages of being analytic, continuous, verifiable, and rapidly convergent, and it can reduce the computation time. The model problem does not need to be linearized or discretized, which is one of the main advantages of this approach. Concerning this, Wazwaz [13], [14] provides the solution of the modified KdV equation and KdV equation which was originally derived by Korteweg and De Vries [2] in the form of a truncated series solution. A generalized KdV equation is numerically solved using the Galerkin method and a quadratic B-spline FEM coupled with ADM by Geyikli and Kaya [15]. Bakodah [16] suggested a new method for solving the generalized fifth-order KdV equation called a modified ADM. The exact solutions and the standard ADM are compared with their numerical solutions, showing that their method works well. An effective ADM-based modification to the convergence parameter for KdV equations was proposed and studied by Bakodah et al. [17].

These conventional methods demand strict adherence to step size requirements and repeated iterations for achieving numerical precision, resulting in substantial computational costs. By integrating ANN models into the study of dynamic problems, researchers aim to overcome these limitations [18]. ANN models offer the potential to approximate solutions with improved efficiency and reduced computational burden. This approach opens up new avenues for understanding complex dynamical systems while alleviating the challenges associated with traditional numerical techniques. As a result, ANN models have emerged as promising tools in the ongoing pursuit of accurate and computationally efficient solutions to dynamic problems.

In response to the challenges posed by traditional numerical techniques, researchers have been concurrently developing ANN models to approximate solutions for dynamic problems. In 1943, McCulloch and Pitts [19] formulated the initial concept of the first basic neural network model. In their model, neurons follow the "all or none" concept for prediction. It attracts researchers to explore and contribute to this emerging area, viz. ANN. However, it became popular among mathematicians in the 90's after Lee and Kang [20] developed the Hopfield neural network. This neural model proved to be a significant breakthrough, particularly for solving ordinary differential equations (ODEs). In a pioneer work, Sahoo and Chakraverty [21] explore the capabilities of artificial neural networks (ANNs) for solving different dynamical problems. In recent work, Ullah et.al [22] developed an ANN model with the power of numerically generated data points for the study of nanofluid flow between two circular plates.

Given the aforementioned challenges and potential benefits of using ANN models for approximating solutions to dynamic problems, it becomes evident that there is a strong motivation to introduce effective intelligent computing based on numerical techniques. These new algorithms aim to provide a better understanding of the dynamic behavior exhibited by GKdV equations. By developing efficient ANN algorithms, researchers can potentially improve the accuracy and computational efficiency of studying such systems, leading to valuable insights and applications in various fields.

## 3. Preliminaries

This section discusses the architecture of the multilayer artificial neural network and the governing equation of the GKdV model.

### 3.1. Architecture of Artificial Neural Network

ANN is a branch of Artificial Intelligence (AI) that mimics the human brain's training process to predict patterns from given historical data. Neural networks are processing devices of mathematical algorithms that may be built using computer languages [18].

Various learning procedures and parameters may be required for modeling a neural network. The neural network is made up of layers, and layers are made up of several neurons /nodes. Every $i^{th}$ and $(i+1)^{th}$ layer neuron is interconnected by synaptic weights. Weights are numerical values allocated to each link. Signals are received by the input layer, which is multiplied by the weights, then summed and sent to one or more hidden layers. The activation function processes its net input to shape the output of the model. The input/output relation is given by

$$O_q = \nabla\left(net\,(q)\right),$$

$$net\,(q) = \sum_{i=1}^{n} (\omega_{qi} . o_i) + \tilde{b}_q \,, q = 1,2,3,....N \,,$$

where $O_q$ is the output for perceptron $q$ after passing through an activation function $\nabla$, $o_i$ denotes the given input in the first layer, $\omega_{qi}$ denotes weight from the input unit $i$ to the hidden unit $q$, and $\tilde{b}_q$ is the bias.

ANN analyses data in the same way as the human brain does. It analyses input and output values of the training dataset and updates the values of weights to improve accuracy through fine-tuning.

### 3.2. Geophysical Korteweg-de Vries equation (GKdV)

Several types of KdV equations are commonly used in geophysics, each with slightly different formulations and applications. These includes the Standard KdV equation, the Modified KdV equation, the Higher-order KdV equation, the Dispersion-managed KdV equation, the generalized KdV equation, and the GKdV equation. However, we are focusing on the GKdV equation widely used in geophysics and oceanography to describe the behavior of surface waves in the ocean and other areas of physics, such as plasma physics and nonlinear optics.

The governing nonlinear GKdV equation [23], [24] can be expressed as;

$$\eta_\tau - w\eta_x + \frac{3}{2}\eta\eta_x + \frac{1}{6}\eta_{xxx} = 0, \tag{1}$$

where,

- $\eta$ is the surface elevation of water,
- $w$ is the Coriolis constant (Coriolis effect: The impact of the rotation of the earth on fluid),
- $x$ and $\tau$ represent the spatial and temporal coordinates, respectively.

### 4. Methodology

In this section, a numerical paradigm-based ANN has been considered to find out the solution to the governing equation (Eq. 1) concerning the following initial condition [25], [26];

$$\eta(x,0) = 2(u+w)\sec h^2\left(\sqrt{1.5(u+w)}\ x\right), \tag{2}$$

where $u$ is the wave speed.

### 4.1 Implementation of ADM on the GKdV equation

Eq. (1) can be transformed into an operator form using ADM [27] as follows:

$$L\eta + R\eta + \frac{3}{2}N\eta = 0, \qquad x \in \Re, \tag{3}$$

where the linear differential operators $L$ *and* $R$ are:

$$L = \frac{\partial}{\partial\tau} \quad \text{and}\ R = \frac{1}{6}\frac{\partial^3}{\partial x^3} - \omega\frac{\partial}{\partial x}\ , \tag{4}$$

and the nonlinear operator is $N\eta$, which can be expressed as $N\eta = \eta\eta_x$

According to ADM, the solution of the governing equation may be formulated as an infinite series of the form

$$\eta(x,\tau)=\sum_{n=0}^{\infty}\eta_n(x,\tau), \tag{5}$$

where $\eta_n(x,\tau)$ will be obtained recurrently.

The nonlinear term $N\eta$ can be decomposed into an infinite series of polynomials, which is given by

$$N\eta=\sum_{n=0}^{\infty}A_n\left(\eta_0,\eta_1,\eta_2,\eta_3,...\eta_n\right), \tag{6}$$

where the components $A_n$ are set of Adomian polynomials and can be defined in the following manner:

$$A_n=\frac{1}{n!}\frac{d^n}{d\mu^n}\left[N\left(\sum_{j=0}^{n}\mu^j\eta_j\right)\right]_{\mu=0} \tag{7}$$

Performing the inverse operator $L^{-1}$ on Eq. (3) and adopting the initial condition given in Eq. (2), we get

$$\eta(x,\tau)=\eta(x,0)-L^{-1}\left\{R\eta+\frac{3}{2}N\eta\right\}. \tag{8}$$

Now, substituting Eqns. (6) and (7) into Eq. (8), gives

$$\eta(x,\tau)=\sum_{n=0}^{\infty}\eta_n(x,\tau)=\eta(x,0)-L^{-1}\left(R\sum_{n=0}^{\infty}\eta_n\right)-L^{-1}\left(\frac{3}{2}\sum_{n=0}^{\infty}A_n\right) \tag{9}$$

where $A_n$ are adomian polynomials that serve as non-linear terms.

Here $N\eta=\eta\eta_x$. So $A_n$ is given by:

$$\begin{aligned}
A_0&=\eta_0(\eta_0)_x,\\
A_1&=\eta_1(\eta_0)_x+\eta_0(\eta_1)_x,\\
A_2&=\eta_2(\eta_0)_x+\eta_1(\eta_1)_x+\eta_0(\eta_2)_x,\\
A_3&=\eta_3(\eta_0)_x+\eta_2(\eta_1)_x+\eta_1(\eta_2)_x+\eta_0(\eta_3)_x,
\end{aligned} \tag{10}$$

.

.

.

Setting the aforesaid polynomials in Eq. (9), the components of $\eta_n(x,\tau)$ are as follows

$$
\begin{aligned}
&\eta_0(x,\tau)=\eta(x,0)\\
&\eta_1(x,\tau)=-L^{-1}\left(R\eta_0+\frac{3}{2}A_0\right),\\
&\eta_2(x,\tau)=-L^{-1}\left(R\eta_1+\frac{3}{2}A_1\right),\\
&\eta_3(x,\tau)=-L^{-1}\left(R\eta_2+\frac{3}{2}A_2\right),\\
&\cdot\\
&\cdot\\
&\cdot\\
&\eta_{n+1}(x,\tau)=-L^{-1}\left(R\eta_n+\frac{3}{2}A_n\right),\qquad n\geq 0.
\end{aligned}
\tag{11}
$$

Combining the above approximation, the solution of Eq. (1) may be obtained in the series form as

$$\eta=\eta_0+\eta_1+\eta_2+\ldots \tag{12}$$

### 4.2 Black box modeling for GKdV equation

This section explains the formation of the undertaken ANN technique to solve the GKdV equation. Here, we have considered a multilayer ANN model having $n$ neurons and $k$ input training datasets with batch size $L$. Then, the forward propagation for the ANN model is defined as

$$X\rightarrow\omega^{(1)T}X\rightarrow\nabla\left(\omega^{(1)T}X\right)=H^{(1)}\rightarrow\omega^{(2)T}H^{(1)}\rightarrow\nabla\left(\omega^{(2)T}H^{(1)}\right)=H^{(2)}....=Output,$$

where $X=\begin{bmatrix} x_1^1 & \cdots & x_1^k & \cdots & x_1^L\\ \vdots & \ddots & \vdots & & \vdots\\ x_i^1 & \cdots & x_i^k & \cdots & x_i^L\\ \vdots & & \vdots & \ddots & \vdots\\ x_m^1 & \cdots & x_m^k & \cdots & x_m^L \end{bmatrix}$ is $(m\times L)$ matrix of input neurons,

$\omega^{(i)}$ is randomly generated weights of order $(m\times n)$ for layer $i$,
$H^{(i)}$ is the output of $i^{th}$ layer.

After the formation of NN, we need to tune the adjustable parameters and minimize the objective function by applying a suitable optimizer. Here, we have taken Adam as an optimizer that combines the concepts of momentum and adaptive learning rates. In addition, *L2* regularization has been used, which acts as a penalty to the error function and helps in keeping weights smaller.

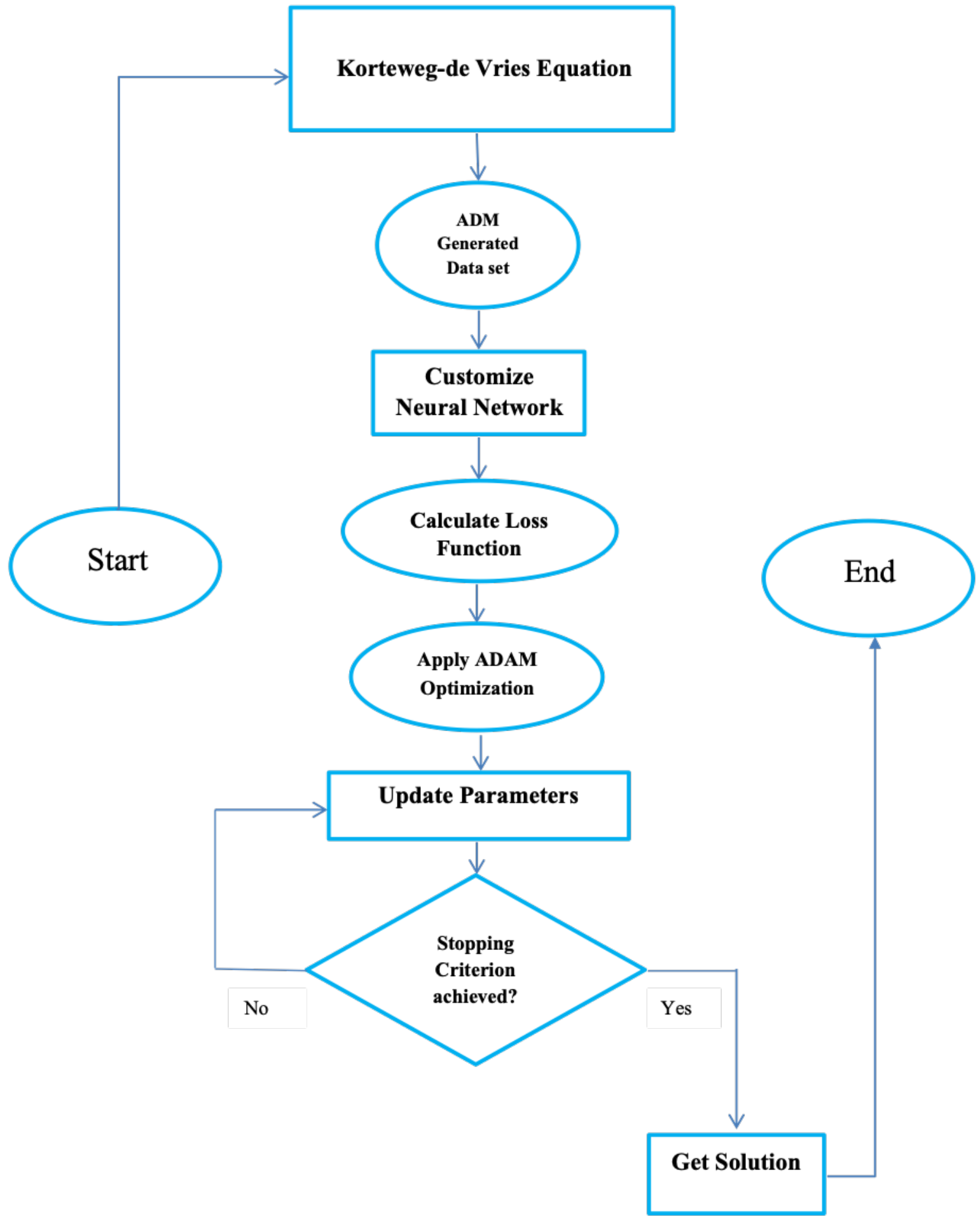

**Figure 1:** Workflow diagram of the G-KdVNet

## 5. Simulation Results and Discussion

This section demonstrated the accuracy and applicability of ADM and ANN techniques via different simulations. In this case, all the calculations have been performed using a finite number of terms $(n=5)$ by truncating the infinite series. We have considered the values of the parameter $u=0.5$ for this calculation.

The exact/analytical solution of the governing GKdV equation (Eq. (1)) is [24]

$$\eta(x,\tau)=2(u+w)\sec h^{2}\left\{\sqrt{1.5(u+w)}\,(\xi-u\tau)\right\} \tag{13}$$

### 5.1 Results Obtained by ADM

In solving the governing equation by using the above-discussed approach, we obtained a series of solutions with up to five terms. In this regard, Figure 2 demonstrates a 3D view of the soliton solutions of the GKdV equation for $w = 0$. One may observe that the solution graph obtained by ADM is quite similar to the exact solution of the GKdV equation given in Eq. (1). In Figure 3, we have also shown the term-by-term convergence 2D plot of the obtained series solution. One may note that there is no difference between the exact solution and the numerical solution made up of four terms. Beyond the first four terms of the series solution, there are no discernible changes in the ADM solution, as seen from the five-term solution. Error plots, which show the absolute difference between approximate and exact solutions, were computed by taking into account different terms of the series solution (Eq. 12). These plots are depicted in Figure 4, and they demonstrate that the first four or five terms may be sufficient to obtain the converged solution for the current problem. Also, it may be noted that as solution terms rise from one to five, the error becomes less.

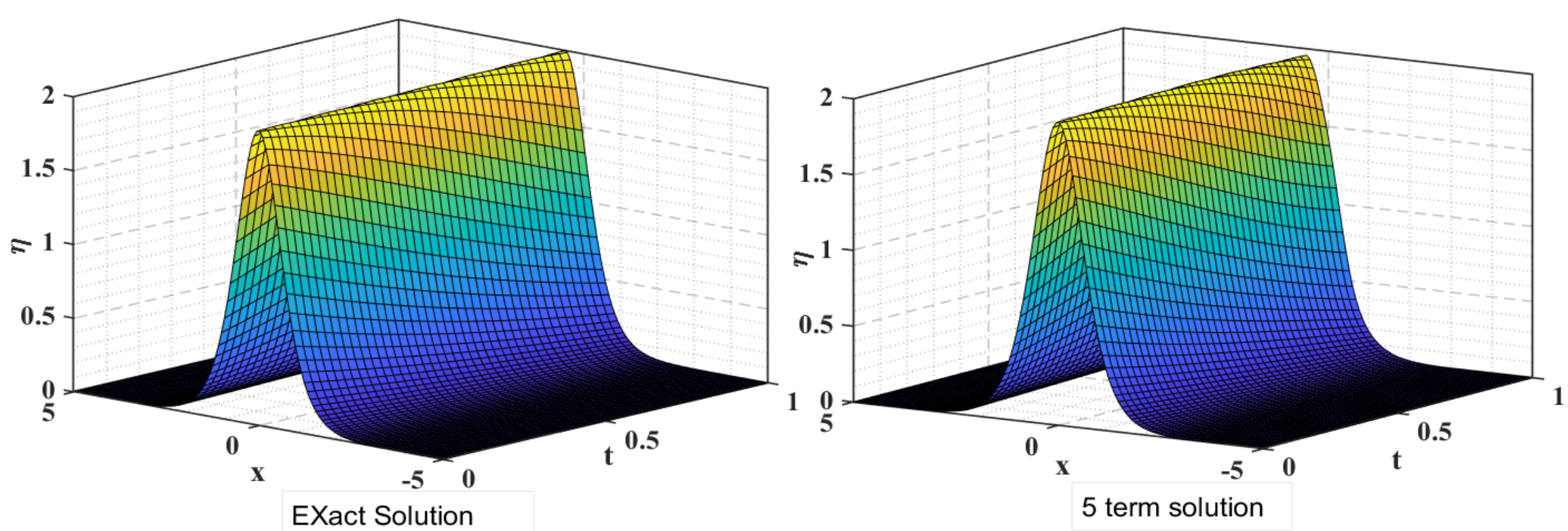

**Figure 2:** Surface plot of $\eta$: Exact solution and Solution by ADM

From Figures 5–6, one may conclude that a significant increase in the Coriolis constant results in a rise in wave height and a shrinking of the wavelength. Similar behavior may be observed in Figure 7, which demonstrates a three-dimensional perspective of the soliton solutions of GKdV equation for the following values of the Coriolis constant: $w = 0.5$, $w = 1$, and $w = 1.5$. The typical changes can be observed from these observations in response to an increase in the value of the Coriolis constant. As a result, we may conclude that the Coriolis constant is inversely proportional to wavelength and directly proportional to wave height. This phenomenon is sometimes referred to as waves having their crests grow more sharpened and their troughs becoming more flattened.

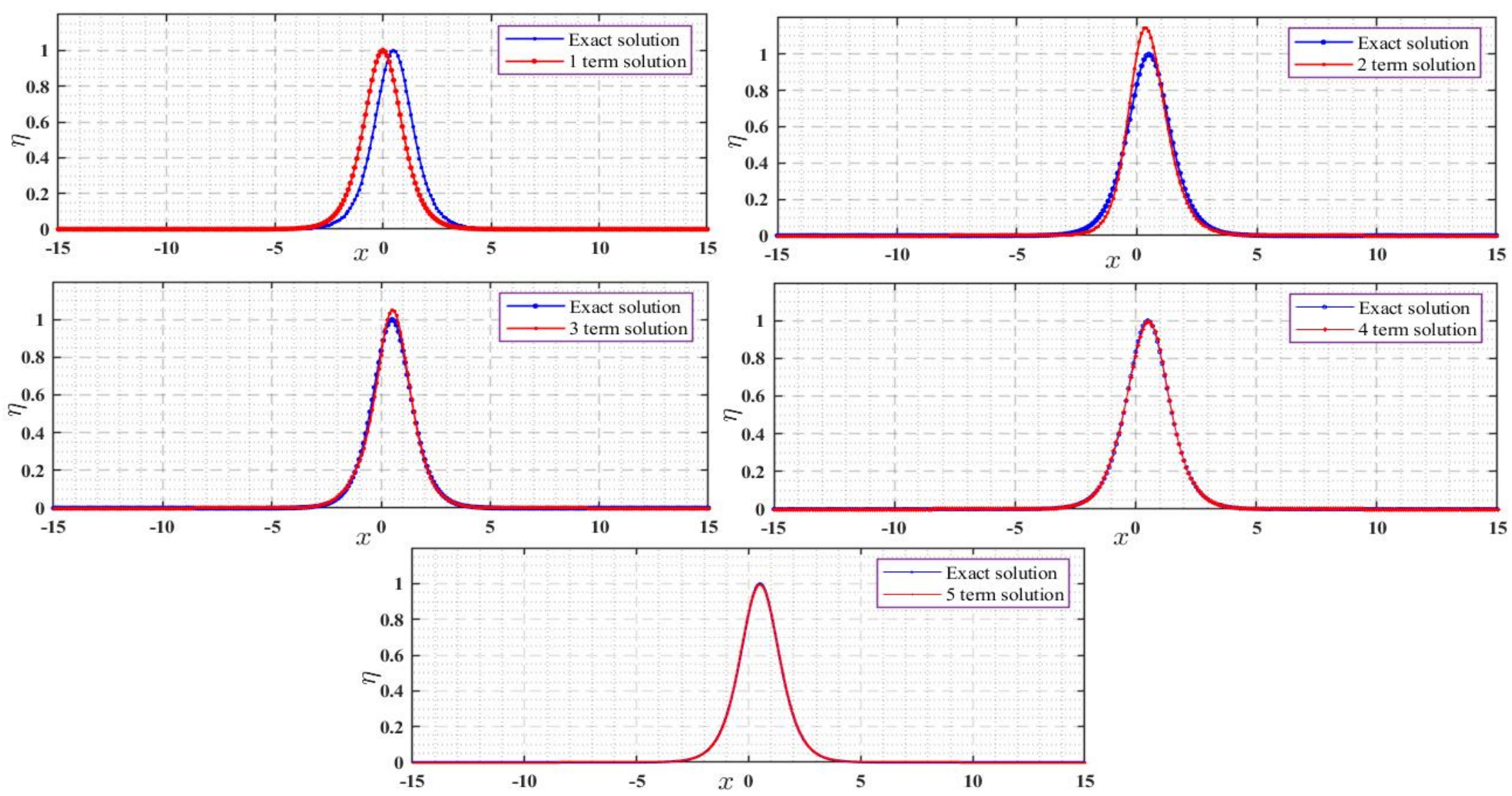

**Figure 3:** Term-by-term convergence plot of $\eta$

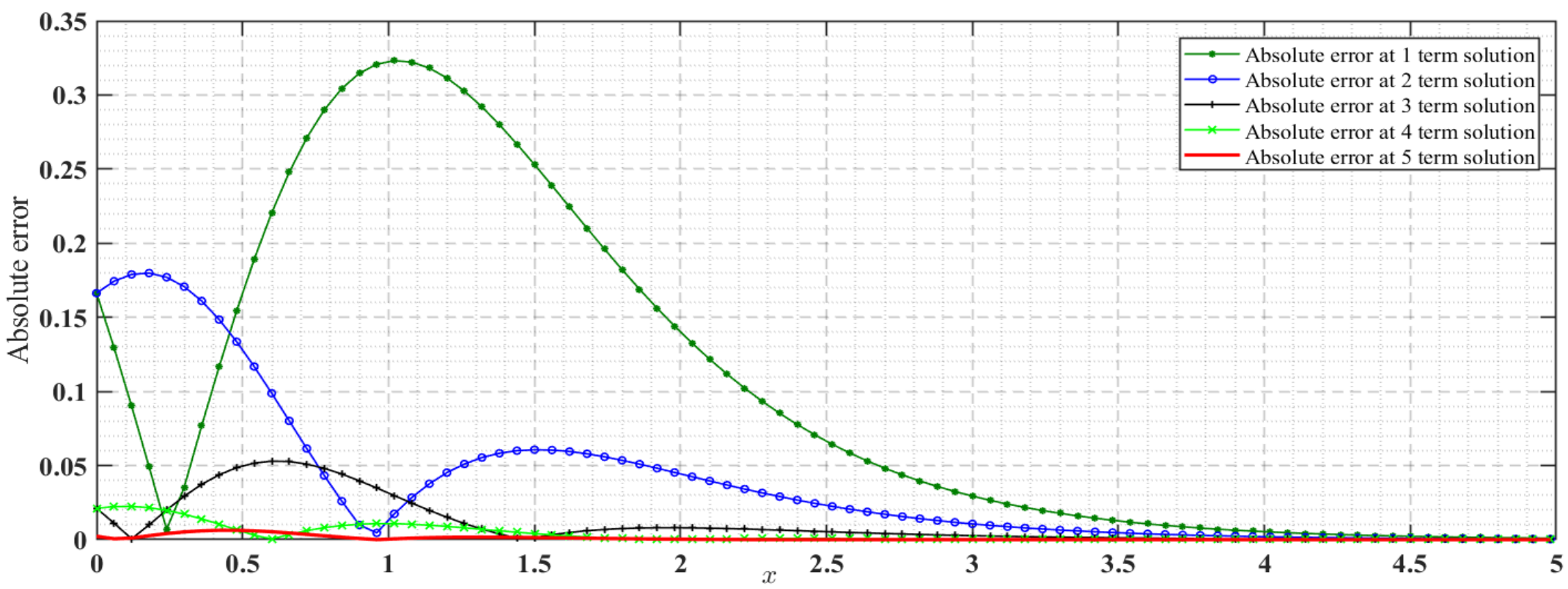

**Figure 4:** Plot of absolute error of GKdV equation by ADM for different terms

The incorporation of the Coriolis term into the GKdV equation results in a discernible change in the shape of the solution. In practice, it follows that the inclusion of the Coriolis constant in the KdV equation has an effect on tsunami wave propagation. [28]

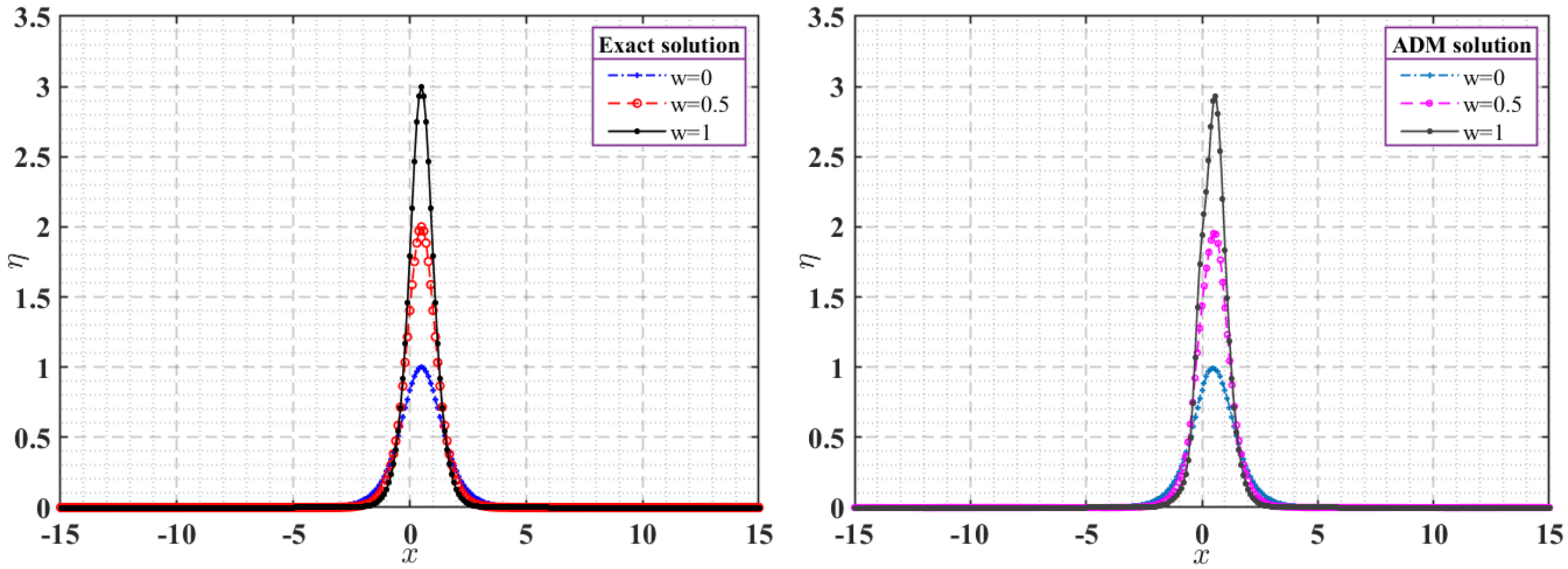

**Figure 5:** Comparison plot of surface elevation $\eta$ of Eq. (1) and Solution by ADM at $w=0, 0.5\, and\, 1,\ u=0.5\ and\ \tau=1$

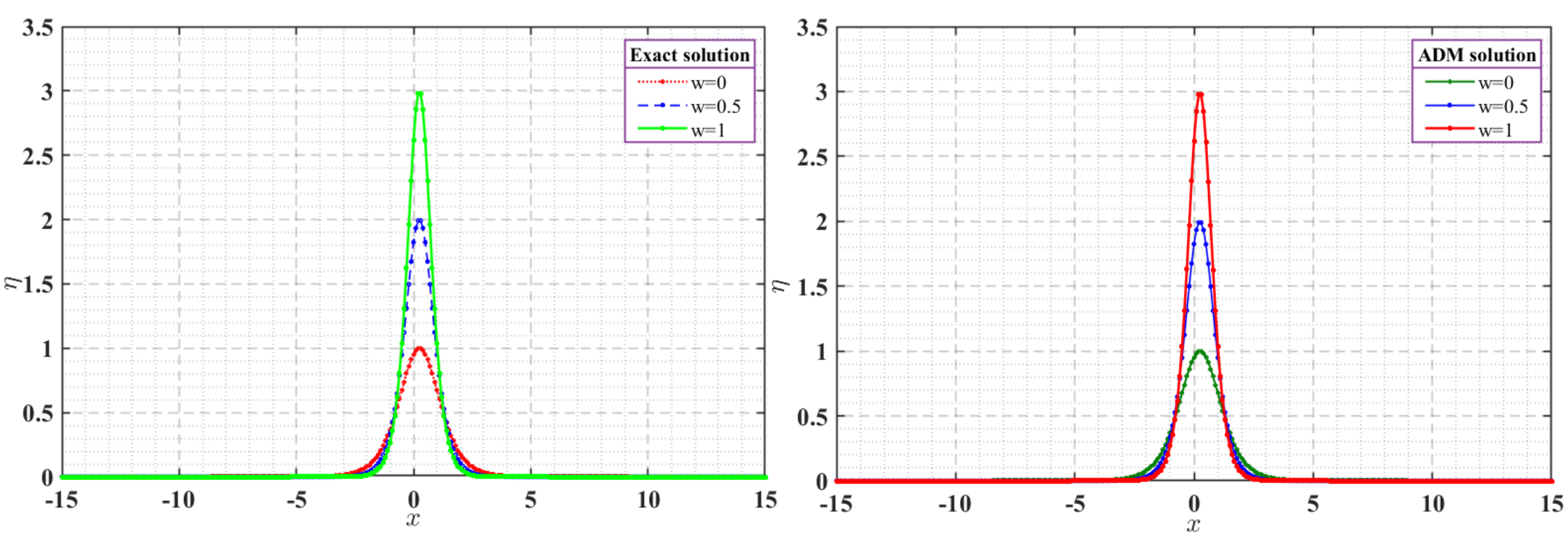

**Figure 6:** Comparison plot of the surface elevation $\eta$ of Eq. (1) and Solution by ADM at $w=0, 0.5\, and\, 1,\ u=0.5\ and\ \tau=0.5$

**Table 1:** Comparison of obtained results vs existing solution [23], when $w=1$

| No. of terms | $x=1$ **and** $\tau=0.1$ | | $x=2$ **and** $\tau=0.5$ | | $x=3$ **and** $\tau=1$ | |
|---|---|---|---|---|---|---|
| | **HPM [23]** | **ADM** | **HPM [23]** | **ADM** | **HPM [23]** | **ADM** |
| 1st | 0.5421 | 0.5421 | 0.0296 | 0.0296 | 0.0015 | 0.0015 |
| 2nd | 0.6157 | 0.6157 | 0.0517 | 0.0517 | 0.0037 | 0.0037 |
| 3rd | 0.6202 | 0.6202 | 0.0599 | 0.0599 | 0.0054 | 0.0054 |
| 4th | 0.6203 | 0.6203 | 0.0619 | 0.0619 | 0.0062 | 0.0062 |
| 5th | 0.6203 | 0.6203 | 0.0623 | 0.0623 | 0.0065 | 0.0065 |
| 6th | 0.6203 | 0.6203 | 0.0623 | 0.0623 | 0.0066 | 0.0066 |

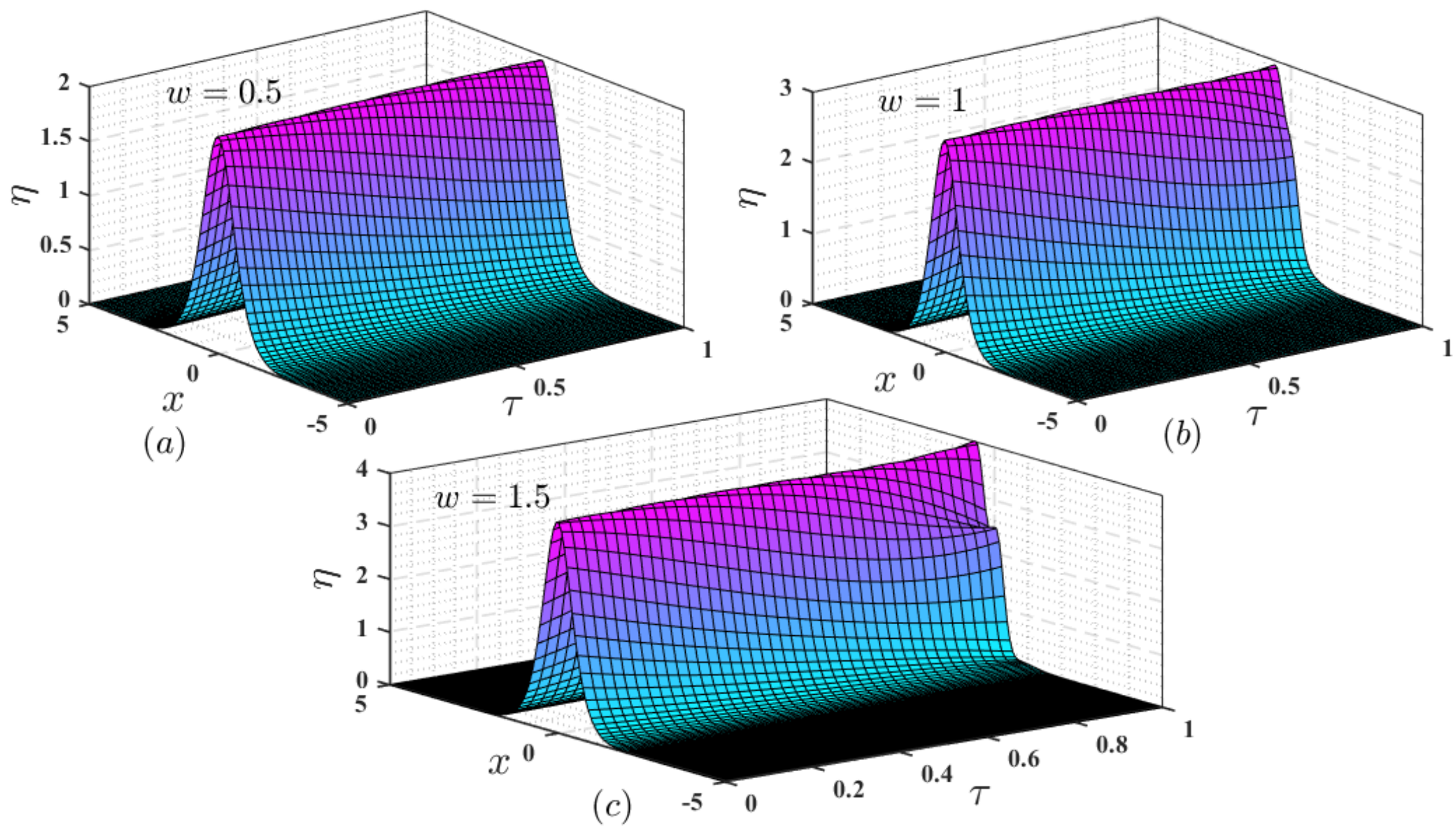

**Figure 7:** Surface Plot of $\eta$ by ADM at **(a)** $w=0.5$, **(b)** $w=1$, **(c)** $w=1.5$.

To validate the obtained solutions, several tabular results are provided. Table 1 presents a comparison between the ADM solution and the solution from [23] using $w=1$ and three different values $x$ and $\tau$, demonstrating good agreement between the two. Additionally, Table 2 extends this comparison by contrasting the exact solution of the GKdV equation with the ADM solution at a fixed $\tau=0.5$, while varying the Coriolis parameter $\left(w=0,0.5\,and\,1\right)$. The results show that while the ADM approach provides reliable solutions, it performs effectively across a range of values for the Coriolis parameter, demonstrating its versatility for different conditions.

Further, absolute and residual errors have been calculated for various values of $x,\tau$ and $w$ to demonstrate the convergence, accuracy, and reliability of the ADM solution (21). Accordingly, Table 3 presents the absolute error for $w=0$, $\tau=1$, and various values $x$, i.e., $x=0,2,4,6,8\;and\;10$. Similarly, Table 4 shows the absolute error for $x=5$, $w=0$, and different values of $\tau$, i.e., $\tau=0,0.2,0.4,0.6,0.8\;and\;1$. Following that, Table 5 provides the residual errors for $x=10$, $w=0$, and several $\tau$ values, i.e., $\tau=0,0.2,0.4,0.6,0.8\;and\;1$. Finally, Table 6 illustrates residual errors for $w=0$ and various $x$ values, i.e., $x=0,2,4,6,8\;and\;10$ at fixed $\tau=1$. These tables clearly show that as the number of terms increases, both absolute and residual errors tend to approach zero, highlighting the high precision and convergence of the ADM solution.

**Table 2:** Exact solution vs. Solution obtained by ADM

| $x$ | $\tau = 0.5$ | | | | | |
|---|---|---|---|---|---|---|
| | $w = 0$ | | $w = 0.5$ | | $w = 1$ | |
| | **Exact** | **ADM** | **Exact** | **ADM** | **Exact** | **ADM** |
| -2.0 | 0.0780 | 0.0780 | 0.0321 | 0.0321 | 0.0140 | 0.0141 |
| -1.8 | 0.1085 | 0.1085 | 0.0521 | 0.0521 | 0.0255 | 0.0256 |
| -1.6 | 0.1499 | 0.1499 | 0.0843 | 0.0843 | 0.0463 | 0.0464 |
| -1.4 | 0.2053 | 0.2053 | 0.1357 | 0.1357 | 0.0838 | 0.0839 |
| -1.2 | 0.2777 | 0.2776 | 0.2168 | 0.2166 | 0.1510 | 0.1509 |
| -1.0 | 0.3693 | 0.3693 | 0.3417 | 0.3412 | 0.2694 | 0.2687 |
| -0.8 | 0.4804 | 0.4804 | 0.5274 | 0.5268 | 0.4728 | 0.4709 |
| -0.6 | 0.6071 | 0.6072 | 0.7885 | 0.7884 | 0.8062 | 0.8038 |
| -0.4 | 0.7398 | 0.7400 | 1.1239 | 1.1253 | 1.3085 | 1.3107 |
| -0.2 | 0.8623 | 0.8625 | 1.4973 | 1.4995 | 1.9619 | 1.9711 |
| 0.0 | 0.9546 | 0.9546 | 1.8236 | 1.8242 | 2.6147 | 2.6177 |
| 0.2 | 0.9981 | 0.9980 | 1.9925 | 1.9906 | 2.9832 | 2.9744 |
| 0.4 | 0.9833 | 0.9831 | 1.9340 | 1.9320 | 2.8531 | 2.8475 |
| 0.6 | 0.9135 | 0.9133 | 1.6732 | 1.6728 | 2.3043 | 2.3055 |
| 0.8 | 0.8035 | 0.8034 | 1.3100 | 1.3106 | 1.6218 | 1.6242 |
| 1.0 | 0.6736 | 0.6736 | 0.9481 | 0.9487 | 1.0351 | 1.0363 |
| 1.2 | 0.5423 | 0.5423 | 0.6480 | 0.6483 | 0.6203 | 0.6206 |
| 1.4 | 0.4226 | 0.4226 | 0.4259 | 0.4260 | 0.3579 | 0.3578 |
| 1.6 | 0.3210 | 0.3211 | 0.2727 | 0.2727 | 0.2020 | 0.2019 |
| 1.8 | 0.2392 | 0.2392 | 0.1718 | 0.1717 | 0.1126 | 0.1125 |
| 2.0 | 0.1757 | 0.1757 | 0.1070 | 0.1070 | 0.0623 | 0.0623 |

**Table 3:** Absolute Error at $\tau = 1$, $w = 0$ and for various $x$ values

| No of terms | $x = 0$ | $x = 2$ | $x = 4$ | $x = 6$ | $x = 8$ | $x = 10$ |
|---|---|---|---|---|---|---|
| 1 | 0.16633 | 0.14014 | 0.00536 | $0.16894 \times 10^{-3}$ | $0.05289 \times 10^{-4}$ | $0.01656 \times 10^{-5}$ |
| 2 | 0.16633 | 0.04438 | 0.00198 | $0.06271 \times 10^{-3}$ | $0.0196 \times 10^{-4}$ | $0.00615 \times 10^{-5}$ |
| 3 | 0.02117 | 0.00803 | 0.00052 | $0.01672 \times 10^{-3}$ | $0.0052 \times 10^{-4}$ | $0.00164 \times 10^{-5}$ |
| 4 | 0.02117 | 0.00029 | 0.00011 | $0.00345 \times 10^{-3}$ | $0.0011 \times 10^{-4}$ | $0.00034 \times 10^{-5}$ |
| 5 | 0.00227 | 0.00032 | 0.00002 | $0.00058 \times 10^{-3}$ | $0.00018 \times 10^{-4}$ | $0.00006 \times 10^{-5}$ |

**Table 4:** Absolute Error at $x = 5$, $w = 0$ and for various $\tau$ values

| No of terms | $\tau = 0$ | $\tau = 0.2$ | $\tau = 0.4$ | $\tau = 0.6$ | $\tau = 0.8$ | $\tau = 1$ |
|---|---|---|---|---|---|---|
| 1 | 0 | $0.09406\times10^{-1}$ | $0.20468\times10^{-1}$ | $0.33453\times10^{-1}$ | $0.48660\times10^{-1}$ | $0.66420\times10^{-1}$ |
| 2 | 0 | $0.07481\times10^{-2}$ | $0.31543\times10^{-2}$ | $0.74844\times10^{-2}$ | $0.14036\times10^{-2}$ | $0.23140\times10^{-2}$ |
| 3 | 0 | $0.00378\times10^{-2}$ | $0.03130\times10^{-2}$ | $0.10915\times10^{-2}$ | $0.26712\times10^{-2}$ | $0.53820\times10^{-2}$ |
| 4 | 0 | $0.00124\times10^{-3}$ | $0.02010\times10^{-3}$ | $0.10286\times10^{-3}$ | $0.32790\times10^{-3}$ | $0.80512\times10^{-3}$ |
| 5 | 0 | $0.00020\times10^{-4}$ | $0.00550\times10^{-4}$ | $0.03909\times10^{-4}$ | $0.151721\times10^{-4}$ | $0.41623\times10^{-4}$ |

**Table 5:** Residual Error at $x = 10$, $w = 0$ and for various $\tau$ values

| No of terms | $\tau = 0.2$ | $\tau = 0.4$ | $\tau = 0.6$ | $\tau = 0.8$ | $\tau = 1$ |
|---|---|---|---|---|---|
| 1 | $-0.10409\times10^{-6}$ | $-0.10409\times10^{-6}$ | $-0.10409\times10^{-6}$ | $-0.10409\times10^{-6}$ | $-0.10409\times10^{-6}$ |
| 2 | $-0.01803\times10^{-6}$ | $-0.03605\times10^{-6}$ | $-0.05408\times10^{-6}$ | $-0.07211\times10^{-6}$ | $-0.09014\times10^{-6}$ |
| 3 | $-0.00156\times10^{-6}$ | $-0.00625\times10^{-6}$ | $-0.01405\times10^{-6}$ | $-0.02498\times10^{-6}$ | $-0.03903\times10^{-6}$ |
| 4 | $-0.00009\times10^{-6}$ | $-0.00072\times10^{-6}$ | $-0.00243\times10^{-6}$ | $-0.00577\times10^{-6}$ | $-0.01127\times10^{-6}$ |
| 5 | $-0.00000\times10^{-6}$ | $-0.00006\times10^{-6}$ | $-0.00032\times10^{-6}$ | $-0.00010\times10^{-6}$ | $-0.00243\times10^{-6}$ |

**Table 6:** Residual Error at $\tau = 1$, $w = 0$ and for various $x$ values

| No of terms | $x = 0$ | $x = 2$ | $x = 4$ | $x = 6$ | $x = 8$ | $x = 10$ |
|---|---|---|---|---|---|---|
| 1 | 0 | - 0.09576 | -0.00338 | $-0.01062\times10^{-2}$ | $-0.03325\times10^{-4}$ | $-0.10409\times10^{-6}$ |
| 2 | 0.375 | - 0.09358 | - 0.00294 | $-0.00920\times10^{-2}$ | $-0.02880\times10^{-4}$ | $-0.09014\times10^{-6}$ |
| 3 | - 0.21094 | -0.04036 | - 0.00128 | $-0.00398\times10^{-2}$ | $-0.01247\times10^{-4}$ | $-0.03903\times10^{-6}$ |
| 4 | - 0.04101 | - 0.00822 | - 0.00037 | $-0.00115\times10^{-2}$ | $-0.00360\times10^{-4}$ | $-0.01127\times10^{-6}$ |
| 5 | 0.07251 | - 0.00041 | - 0.00008 | $-0.00025\times10^{-2}$ | $-0.00078\times10^{-4}$ | $-0.00244\times10^{-6}$ |

### 5.2 Results Obtained by G-KdVNet

In this section, G-KdVNet model has been implemented, and its accuracy has been evaluated using various error metrics.

To evaluate the accuracy of the model, we compared the results obtained from the ANN algorithm with traditional numerical results using various statistical measures. The implementation of ANN structure was carried out in the Jupyter Notebook environment using Python 3.0. For selecting an optimal number of hidden layers and the number of nodes in each hidden layer fine tuning has been done. The findings are supported by illustrative 2D and 3D graphs.

For the simulation, a fully connected multi-layer neural network was constructed with three input nodes and a single output node. The network comprises three hidden layers each having 32 neurons. The Adam optimizer with a learning rate of 0.005 was utilized to update the parameters in the test problems. The training dataset contains 183 points in the temporal-spatial domain for different values of the Coriolis constant ranging from $w=0$ to $w=1$. During the training, the learning performance of each epoch was monitored by tracking the training loss. Figure 8 displays the training loss over 1990 epochs, with the last 100 epochs depicted in the subfigures. The training loss consistently decreases throughout the learning process, indicating the robustness of the model. This observation demonstrates that the neural network effectively learns from the data and generalizes well to the validation set. The decreasing trend in the loss functions reflects the model's ability to accurately predict the output for the given input data. To validate the model, it was tested on 20 validation points. Errors relative to both the ADM and exact solutions were graphically represented in a 20-bin histogram, as shown in Figure 9. This demonstrates the promising potential of the constructed neural network for solving the GKdV model.

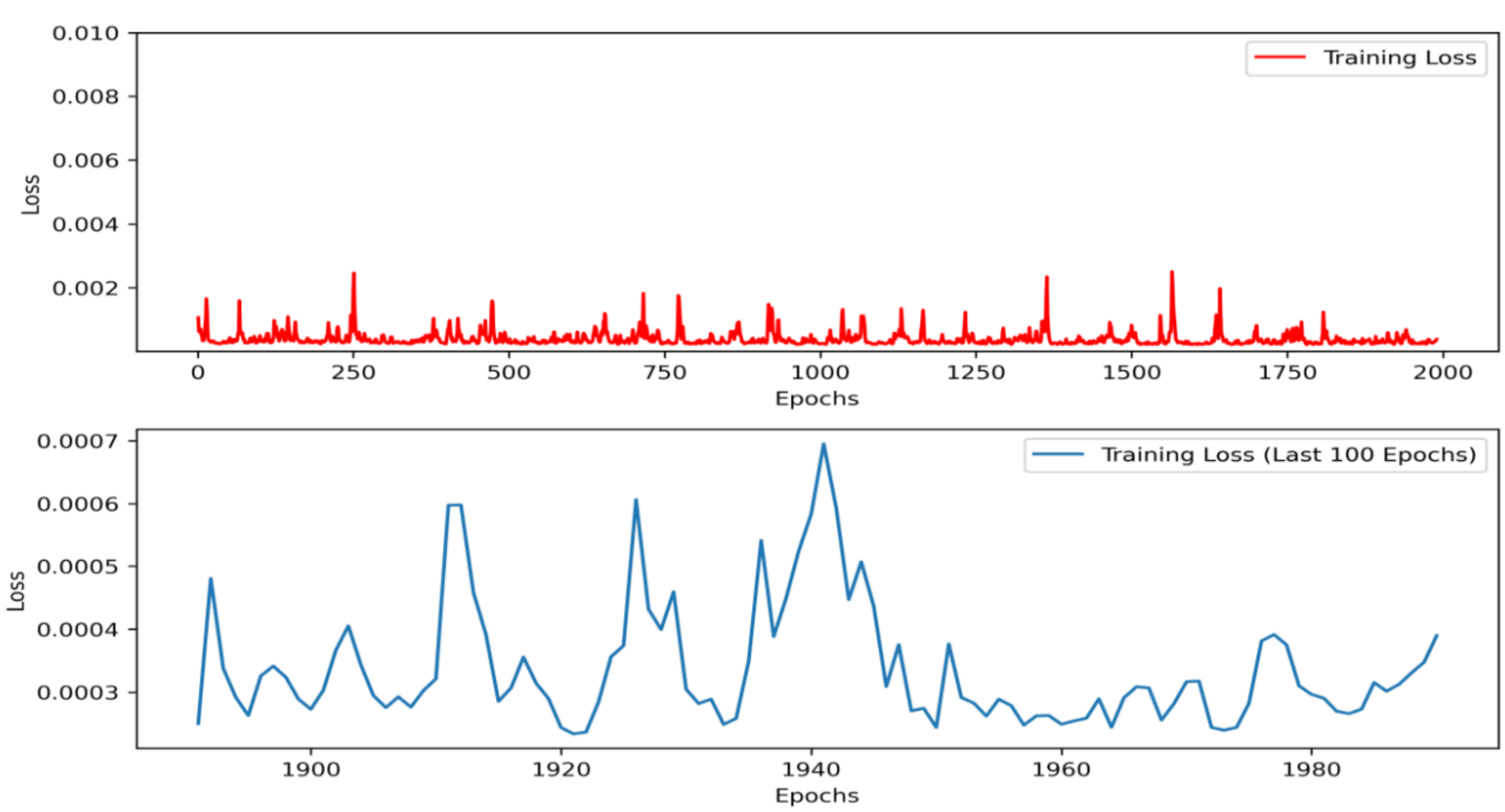

**Figure 8:** Learning curves representing the training loss of G-KdVNet with respect to 1,990 epochs (upper), zoomed portion for the last 100 epochs (lower) during the training process

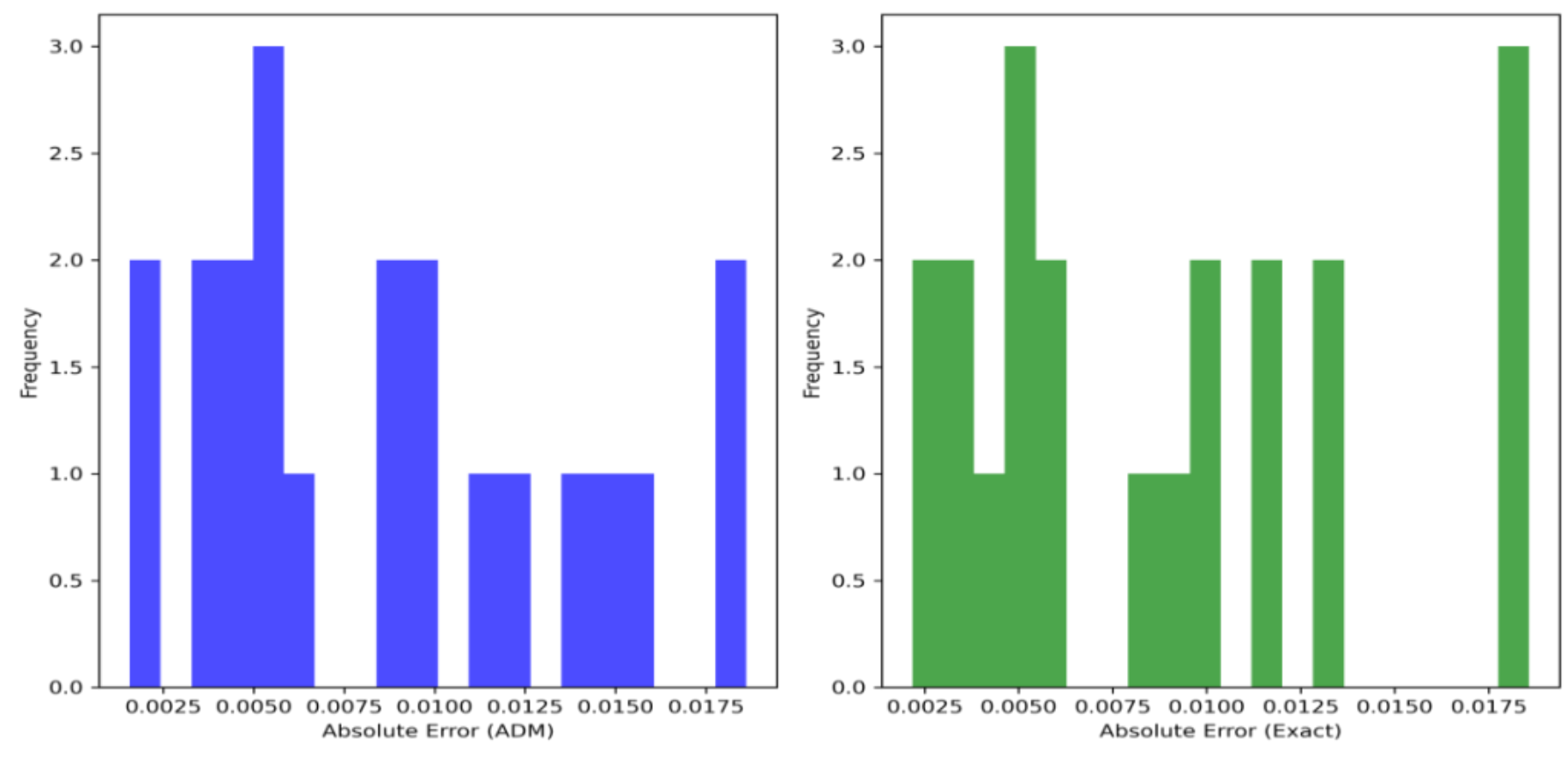

**Figure 9:** 20-bins histogram for absolute errors at 20 validation points

The performance of the ANN model is assessed by comparing its predictions with the ADM and exact (ground truth) solution for the GKdV equation at different time points. Figure 10 represents the heat map comparison of ANN and ADM solutions at the testing dataset with a space-time grid of (21, 6) for different Coriolis constants. The 3D surface plots comparison has been demonstrated in Figure 11, which showcases the influence of various Coriolis constants. These figures likely show how different factors affect the behavior of the system under consideration. From these figures, one may observe that ANN is accurately predicted in the given domain.

For quantitative assessment, the absolute error was measured on the above testing dataset and the mean absolute error has been found to be 1.316201E-04. By leveraging knowledge from numerical simulations, ANN models can expedite the convergence during training and reduce computational time. Nonetheless, the results indicated promising accuracy in ANN's predictions of GKdV equation.

The model has been tested for 10 randomly generated points, and the absolute errors have been presented in 3D and box plot in Figure 12. From this figure, one can observe that at different testing points the errors lie in a neighbourhood of 0. Figure 13 represents a comparison of absolute errors with respect to ADM and exact solutions. The overlapping of errors between the reference and obtained solutions manifests the efficiency of our model.

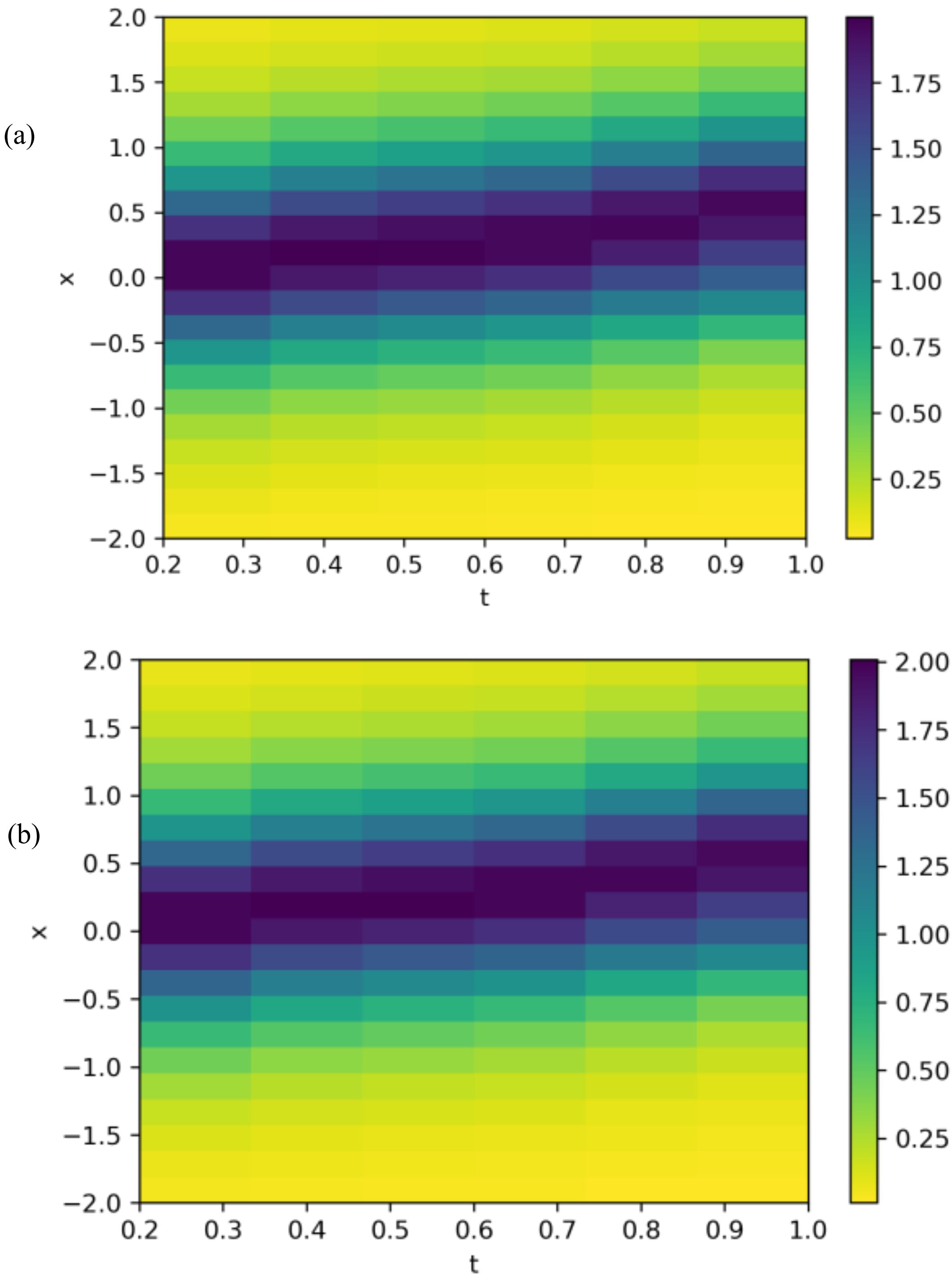

**Figure 10:** 2D plots for comparison of **(a)** ADM and **(b)** G-KdVNet solutions

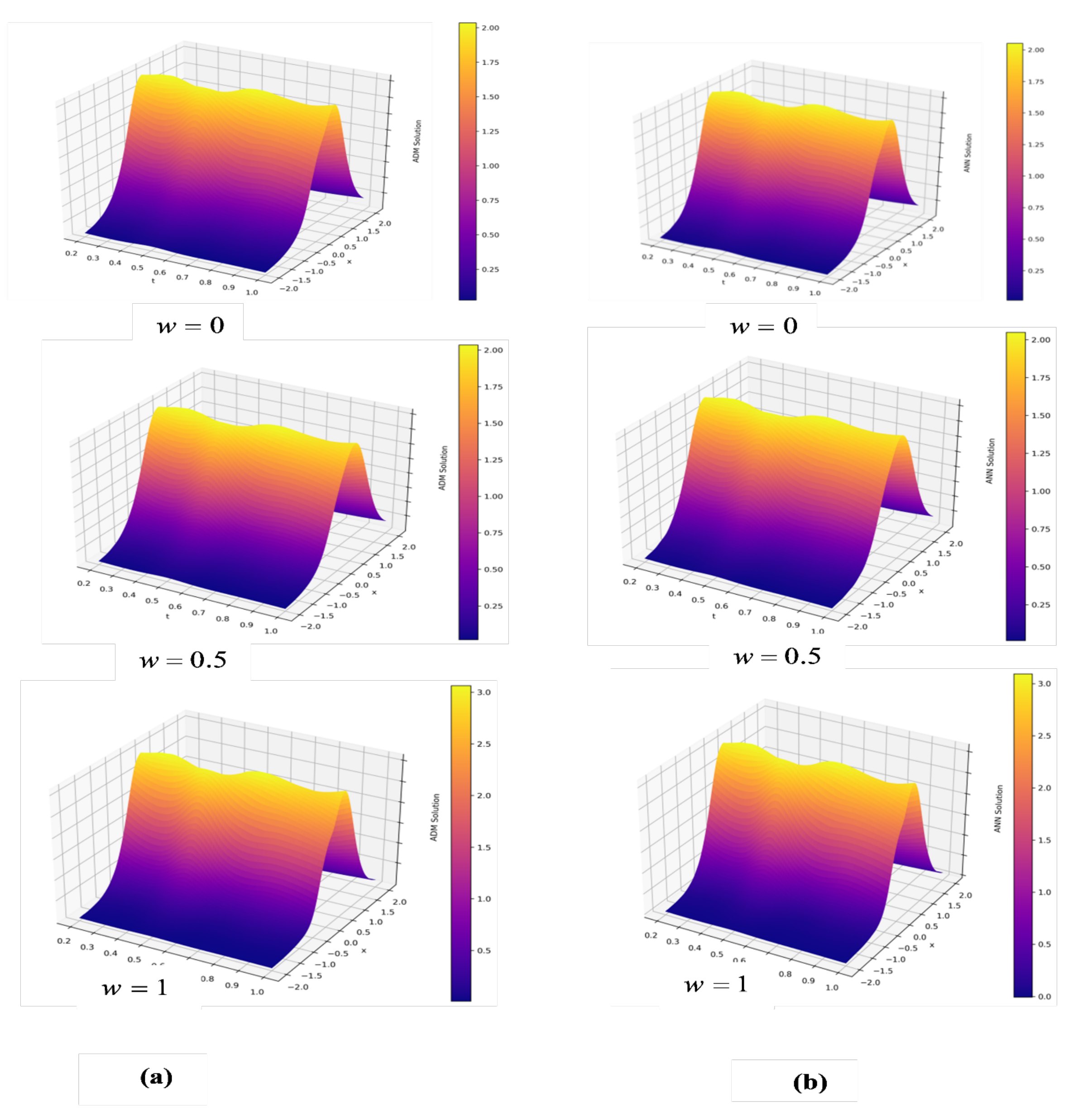

**Figure 11:** 3D surface plots for comparison of **(a)** ADM and **(b)** G-KdVNet solutions for $w = 0, 0.5$ *and* $1$

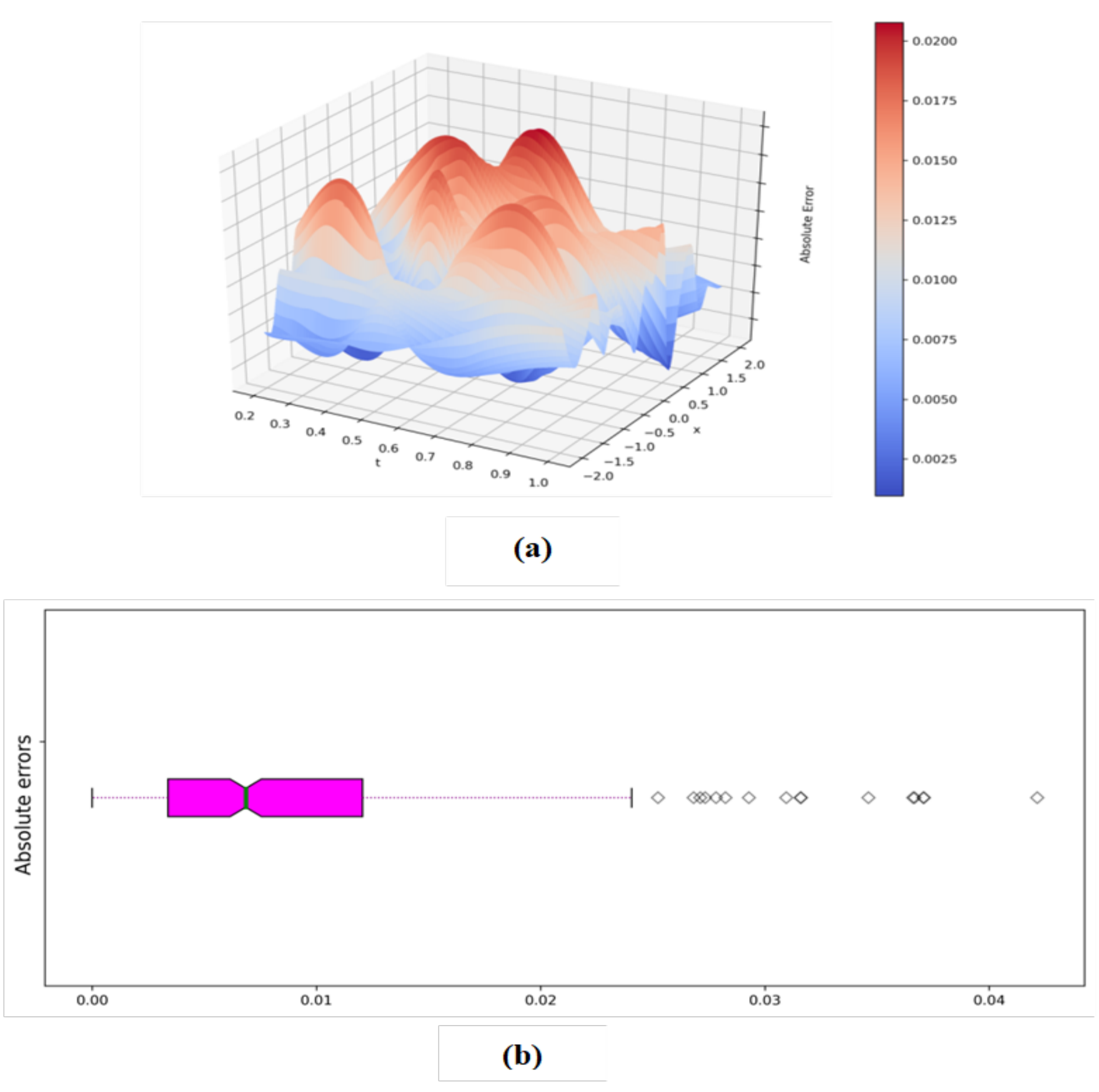

**Figure 12 (a)** Absolute Errors plots in 3D surface, **(b)** Box plot of Absolute Errors between ADM and G-KdVNet solutions at different testing points in the domain

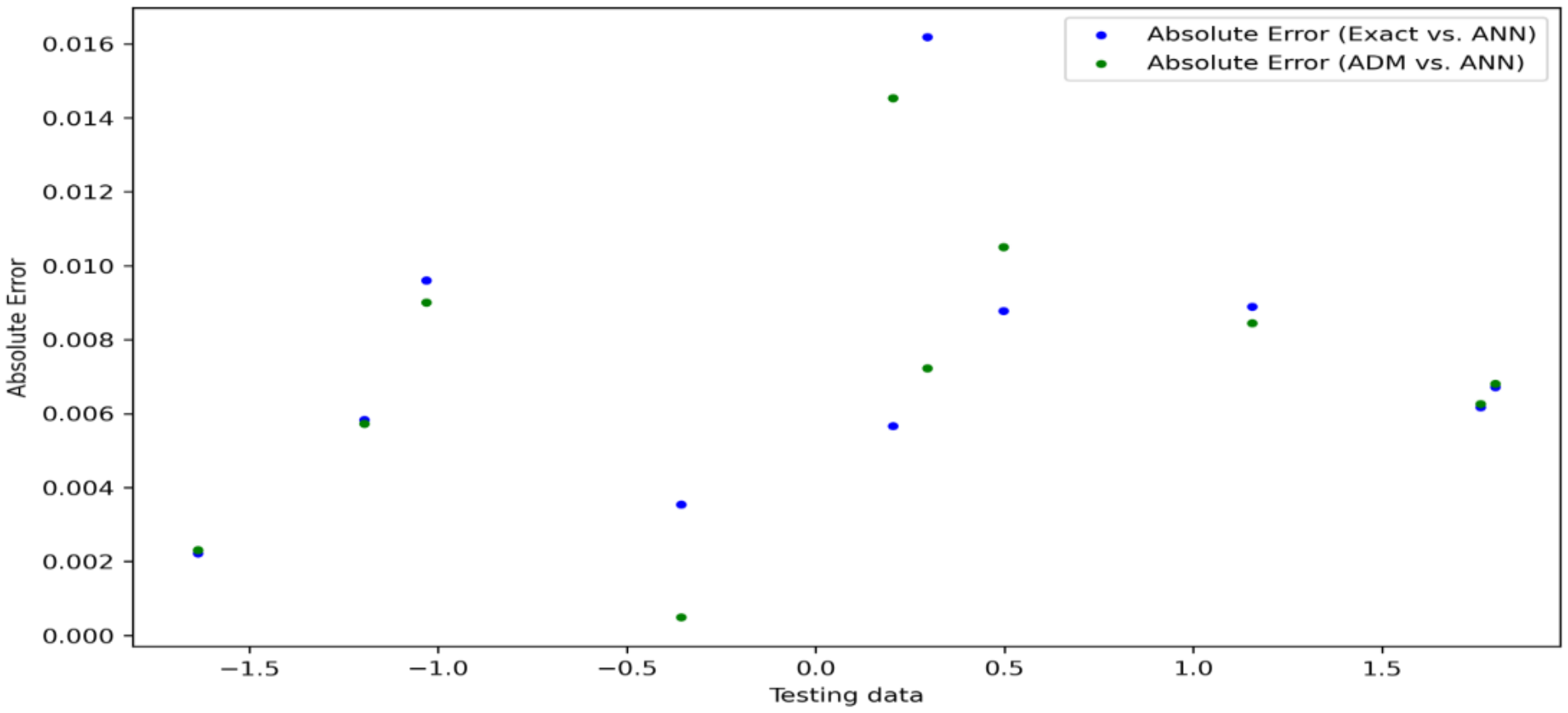

**Figure 13:** Comparison of absolute errors at 10 randomly generated testing points at $\tau = 0.5$, $w = 1$ and for various $x$ values

## 6. Conclusion

The present work investigates the application of effective intelligent computing based on a numerical technique to approximate the GKdV equation that depends upon the Coriolis constant. The obtained results have been validated by comparing them with the exact solution. A key limitation of the ADM method is its lack of precision, particularly over the global solution interval. To address this, a hybrid approach viz. G-KdVNet combining ADM with an ANN model has been implemented, leading to a significant enhancement in the accuracy of the results.

Additionally, a direct relationship between the Coriolis constant and wave height as well as an inverse relationship between the Coriolis constant and wavelength are also observed. The numerical and graphical illustrations provide strong evidence of the effectiveness and convergence of the ANN approach for solving the GKdV equation with various parameters. The ANN approach allows for a comprehensive understanding of the system's response to changes in the physical parameters, providing insights into its dynamics. Overall, the combination of numerical and graphical illustrations contributes to the robustness of the G-KdVNet method and enhances our understanding of the complex system represented by ANN. Nonetheless, the results indicated promising accuracy in G-KdVNet's predictions.

**Authorship Contribution Statement**

M. S: Writing original draft, Methodology.
A. S: Editing, Supervision, and Software.
S. C: Methodology, Writing.

**Acknowledgement**

M.S. would like to express his gratitude to the University Grants Commission (UGC), New Delhi, India, for the support provided in conducting this research.

**Funding**

No funding has received for the preparation of this manuscript.

**Data Availability**

No data has been used for the research described in the article.

**Declaration**

The authors have no conflicts of interest to declare.

**References:**

[1] A. Sahoo, S. Kumar, and S. Chakraverty, "Physics-informed neural network for vibration analysis of large membranes." *Journal of Nonlinear, Complex and Data Science, vol.* 25, no. 7-8, pp. 505-521, 2025.

[2] D. J. Korteweg and G. de Vries, "XLI. *On the change of form of long waves advancing in a rectangular canal, and on a new type of long stationary waves*," *The London, Edinburgh, and Dublin Philosophical Magazine and Journal of Science*, vol. 39, no. 240, pp. 422–443, May 1895, doi: 10.1080/14786449508620739.

[3] R. S. JOHNSON, "Camassa–Holm, Korteweg–de Vries and related models for water waves," *J Fluid Mech*, vol. 455, pp. 63–82, Mar. 2002, doi: 10.1017/S0022112001007224.

[4] N. A. Kudryashov, "On 'new travelling wave solutions' of the KdV and the KdV–Burgers equations," *Commun Nonlinear Sci Numer Simul*, vol. 14, no. 5, pp. 1891–1900, May 2009, doi: 10.1016/J.CNSNS.2008.09.020.

[5] A. R. Seadawy, D. Lu, and C. Yue, "Travelling wave solutions of the generalized nonlinear fifth-order KdV water wave equations and its stability," *Journal of Taibah University for Science*, vol. 11, no. 4, pp. 623–633, Jul. 2017, doi: 10.1016/j.jtusci.2016.06.002.

[6] J. Cai, C. Bai, and H. Zhang, "Efficient schemes for the coupled Schrödinger–KdV equations: Decoupled and conserving three invariants," *Appl Math Lett*, vol. 86, pp. 200–207, Dec. 2018, doi: 10.1016/J.AML.2018.06.038.

[7] P. Karunakar and S. Chakraverty, "Solutions of time-fractional third- and fifth-order Korteweg–de-Vries equations using homotopy perturbation transform method," *Eng Comput (Swansea)*, vol. ahead-of-print, no. ahead-of-print, Jul. 2019, doi: 10.1108/EC-01-2019-0012.

[8] N. H. Aljahdaly, A. R. Seadawy, and W. A. Albarakati, "Analytical wave solution for the generalized nonlinear seventh-order KdV dynamical equations arising in shallow water waves," *Modern Physics Letters B*, vol. 34, no. 26, p. 2050279, Sep. 2020, doi: 10.1142/S0217984920502796.

[9] M. Sahoo and S. Chakraverty, "Sawi Transform Based Homotopy Perturbation Method for Solving Shallow Water Wave Equations in Fuzzy Environment," *Mathematics*, vol. 10, no. 16, p. 2900, Aug. 2022, doi: 10.3390/math10162900.

[10] M. Sahoo and S. Chakraverty, "Dynamics of tsunami wave propagation in uncertain environment," *Computational and Applied Mathematics*, vol. 43, no. 5, p. 266, Jul. 2024, doi: 10.1007/s40314-024-02776-6.

[11] G. Adomian, "A review of the decomposition method and some recent results for nonlinear equations," *Math Comput Model*, vol. 13, no. 7, pp. 17–43, Jan. 1990, doi: 10.1016/0895-7177(90)90125-7.

[12] G. Adomian, *Solving frontier problems of physics: the decomposition method*, vol. 60. Springer Science & Business Media., 2013.

[13] A. M. Wazwaz, "Solitary wave solutions for the modified KdV equation by Adomian decomposition method.," *Int J Appl Math (Sofia)*, vol. 3, no. 4, pp. 361–368, 2000.

[14] A. M. Wazwaz, "Construction of solitary wave solutions and rational solutions for the KdV equation by Adomian decomposition method," *Chaos Solitons Fractals*, vol. 12, no. 12, pp. 2283–2293, Sep. 2001, doi: 10.1016/S0960-0779(00)00188-0.

[15] T. Geyikli and D. Kaya, "Comparison of the solutions obtained by B-spline FEM and ADM of KdV equation," *Appl Math Comput*, vol. 169, no. 1, pp. 146–156, Oct. 2005, doi: 10.1016/J.AMC.2004.10.045.

[16] H. O. Bakodah, "Modified Adomain Decomposition Method for the Generalized Fifth Order KdV Equations," *American Journal of Computational Mathematics*, vol. 03, no. 01, pp. 53–58, 2013, doi: 10.4236/ajcm.2013.31008.

[17] H. O. Bakodah, M. A. Banaja, B. A. Alrigi, A. Ebaid, and R. Rach, "An efficient modification of the decomposition method with a convergence parameter for solving Korteweg de Vries equations," *J King Saud Univ Sci*, vol. 31, no. 4, pp. 1424–1430, Oct. 2019, doi: 10.1016/J.JKSUS.2018.11.010.

[18] A. K. Sahoo and S. Chakraverty, "Machine intelligence in dynamical systems: \A state-of-art review," *WIREs Data Mining and Knowledge Discovery*, vol. 12, no. 4, Jul. 2022, doi: 10.1002/widm.1461.

[19] W. S. McCulloch and W. Pitts, "A logical calculus of the ideas immanent in nervous activity," *Bull Math Biol*, vol. 52, no. 1–2, pp. 99–115, Jan. 1990, doi: 10.1007/BF02459570.

[20] H. Lee and I. S. Kang, "Neural algorithm for solving differential equations," *J Comput Phys*, vol. 91, no. 1, pp. 110–131, Nov. 1990, doi: 10.1016/0021-9991(90)90007-N.

[21] A. K. Sahoo and S. Chakraverty, "A neural network approach for the solution of Van der Pol-Mathieu-Duffing oscillator model," *Evol Intell*, Feb. 2023, doi: 10.1007/s12065-023-00835-1.

[22] H. Ullah *et al.*, "Intelligent computing paradigm for Second-Grade Fluid in a Rotating Frame in a Fractal Porous Medium," *Fractals*, Jul. 2023, doi: 10.1142/S0218348X23401758.

[23] P. Karunakar and S. Chakraverty, "Effect of Coriolis constant on Geophysical Korteweg-de Vries equation," *Journal of Ocean Engineering and Science*, vol. 4, no. 2, pp. 113–121, Jun. 2019, doi: 10.1016/J.JOES.2019.02.002.

[24] A. Geyer and R. Quirchmayr, "Shallow water equations for equatorial tsunami waves," *Philosophical Transactions of the Royal Society A: Mathematical, Physical and Engineering Sciences*, vol. 376, no. 2111, p. 20170100, Jan. 2018, doi: 10.1098/rsta.2017.0100.

[25] R. S. Johnson, *A modern introduction to the mathematical theory of water waves*, vol. 19. Cambridge university press, 1997.

[26] A.-M. Wazwaz, “Negative-order KdV equations in (3+1) dimensions by using the KdV recursion operator,” *Waves in Random and Complex Media*, vol. 27, no. 4, pp. 768–778, Oct. 2017, doi: 10.1080/17455030.2017.1317115.

[27] G. Adomian, “A review of the decomposition method and some recent results for nonlinear equations,” *Math Comput Model*, vol. 13, no. 7, pp. 17–43, Jan. 1990, doi: 10.1016/0895-7177(90)90125-7.

[28] T. Ak, A. Saha, S. Dhawan, and A. H. Kara, “Investigation of Coriolis effect on oceanic flows and its bifurcation via geophysical Korteweg–de Vries equation,” *Numer Methods Partial Differ Equ*, vol. 36, no. 6, pp. 1234–1253, Nov. 2020, doi: 10.1002/num.22469.